\documentclass[a4paper,11pt]{article}
\usepackage[top=30truemm,bottom=30truemm,left=30truemm,right=30truemm]{geometry}
\usepackage{amsmath}
\usepackage{amssymb}
\usepackage{float}
\usepackage{color}
\usepackage[dvipdfmx]{graphicx}
\usepackage{graphics}
\usepackage{graphicx}
\usepackage{tcolorbox}
\usepackage{latexsym}
\usepackage{ascmac}
\usepackage{framed}
\usepackage{fancybox}
\usepackage{theorem}
\usepackage{stmaryrd}
\usepackage{here}
\usepackage{mathtools}
\usepackage{bm}
\usepackage{appendix}

\mathtoolsset{showonlyrefs=true} %%%%%%%%%% 相互参照した式番号のみ表示するためのおまじない

\theoremstyle{break}
\def\qed{\hfill$\Box$}

\newtheorem{defn}{Definition}

\newtheorem{thm}{Theorem}

\makeatletter
\@addtoreset{equation}{section}

\makeatother

\title{{\bf{
Minimal speed of unbounded traveling wave solutions for a 1D reaction-diffusion equation and their relationship with the dynamics at infinity
}}
}

\author{
Yu Ichida
\thanks{
Department of Mathematics, School of Science and Technology, Meiji University, 1-1-1 Higashimita, Tama-ku, Kawasaki 214-8571, Japan, {\tt ichidayu@meiji.ac.jp} }
}

\begin{document}
\maketitle

\begin{abstract}
This paper presents results on the unboundedness and minimal speed of traveling wave solutions for a one-dimensional spatial reaction-diffusion equation with an asymptotically linear reaction term and a saturation parameter.
By applying a Poincar\'e-type compactification, we reveal the full dynamics (including infinity) of the two-dimensional system of ordinary differential equations satisfied by traveling wave solutions.
This yields essential information characterizing traveling wave solutions: the classification of trajectories in the phase plane, the positivity and unboundedness of front-type and sign-changing profiles, and the explicit form of the minimal speed.
This paper examines a special equation with an asymptotically linear reaction term.
While, our results differ from those of conventional linear determinacy. 
We claim that the minimal speed is derived from information at infinity within the traveling wave system.
\end{abstract}


\begin{center}
\end{center}

\section{Introduction}
\label{sec:UBTW-int}
%方程式と記号の説明
In this paper, we are concerned with the following one-dimensional spatial parabolic equation:
\begin{equation}
u_{t}=u_{xx}+\dfrac{u^{3}}{1+su^{2}}, \quad
t>0, \quad x\in \mathbb{R}.
\label{eq:UBTW-int1}
\end{equation}
Let $u = u(t,x)$, $u_{t}=\partial u/\partial t$, and $u_{xx}=\partial^{2} u/\partial x^{2}$.
The parameter $s > 0$ is a constant.

%方程式の背景：物理的背景
According to \cite{FerMai}, the equation \eqref{eq:UBTW-int1} is derived from the propagation of pulses and solitons in optical fibers.
The constant $s>0$ corresponds to the saturation parameter of this phenomenon.
For more details on the physical background, see \cite{FerMai} and the references therein.
%先行研究：grow-upの先行研究紹介
They consider \eqref{eq:UBTW-int1} on an open, smooth, bounded domain in $N$-dimensional space ($N \geq 2$).
Under a setting that generalizes the characteristics of the reaction term, which possesses asymptotic linearity, they prove that a grow-up corresponding to a blow-up at infinite time occurs.

%考える問題について：進行波，非有界進行波，最小速度，線形予測
In this paper, we study the one-dimensional version of this equation on entire space and examine a traveling wave solution as a specific classes of solutions,  considering the equation's physical context.
Our main results show that there are no front-type traveling waves for this equation and that existing traveling waves are unbounded.
Specifically, we provide results regarding the existence and profiles of positive, unbounded traveling wave solutions to the equation \eqref{eq:UBTW-int1}, including their minimal speed. 
We also provide results concerning unbounded traveling wave solutions with sign changes.
Recently, Ito-Ninomiya (\cite{ito1}) provided general and versatile results regarding the minimal speed formula for unbounded traveling wave solutions of reaction-diffusion equations.
However, as will be discussed later, the equation \eqref{eq:UBTW-int1} considered in this paper does not fit within this framework.
In this paper, we present an analysis of the minimal speed and linear determinacy (linear prediction) of unbounded traveling wave solutions in specific cases by investigating the dynamics at infinity of the corresponding second-order ordinary differential equation satisfied by traveling wave coordinates. 
We hope this analysis will serve as a motivation for further theoretical research on unbounded traveling wave solutions.

A traveling wave solution is a typical and representative particular solution to a partial differential equation of the form 
\begin{equation}
u(t,x)=\phi(\xi), \quad \xi=x-ct, \quad c>0,
\label{eq:UBTW-int2}
\end{equation}
where $c>0$ is the wave speed.
The function $\phi$ is a profile, which is a solution that propagates at a constant speed without changing its shape.
Since it is a dynamic and characteristic solution that appears in the spatiotemporal dynamics of solutions to reaction-diffusion equations such as \eqref{eq:UBTW-int1}, many researchers have previously established numerous results regarding its existence, minimal speed, shape, stability, and other properties for a wide range of reaction-diffusion equations.
From \eqref{eq:UBTW-int1}, $\phi(\xi)$ satisfies the following:
\begin{equation}
-c\phi' =\phi''+\dfrac{\phi^{3}}{1+s\phi^{2}}
\label{eq:UBTW-int3}
\end{equation}
with $'=d/d\xi$ and $''=d^{2}/d\xi^{2}$.

An unbounded traveling wave solution is a function $\phi(\xi)$ defined on $\xi \in\mathbb{R}$ that satisfies \eqref{eq:UBTW-int3} and, for some $\xi_{*}\in \mathbb{R}$, satisfies 
\begin{equation}
\lim_{\xi \to \xi_{*}} |\phi(\xi)| =+\infty.
\label{eq:UBTW-int4}
\end{equation}
The minimal speed of an unbounded traveling wave solution is a condition for the existence of a nonnegative unbounded traveling wave solution.
That is, if we denote the minimal speed by $c_{*}$, then this means that a nonnegative unbounded traveling wave solution exists if $c \ge c_{*}$.

When determining the minimal speed under the assumption that the wave propagates in a direction where $\xi$ increases, a linear conjecture can be derived from the information in the linearized equation at the point where $\xi$ approaches infinity.
In other words, the conjecture is that the minimal speed is determined by the linearized eigenvalue information at the stable equilibria in the phase space induced by \eqref{eq:UBTW-int3}.
Determining the minimal speed of traveling waves has been a central focus of numerous studies on reaction-diffusion equations, as it provides crucial information for characterizing the overall behavior of the solutions.

%進行波や最小速度について先行研究列挙しつつ
Extensive research has been conducted to determine the minimal speed of solutions that satisfy the condition:
\begin{equation}
\phi(-\infty)=1, \quad \phi(+\infty)=0,
\label{eq:UBTW-int5}
\end{equation}
where the profile $\phi$ is a monotonically decreasing front-type traveling wave solution.
As the most typical case, consider the one-dimensional reaction-diffusion equation
\begin{equation}
u_{t}=u_{xx}+f(u), \quad t>0, \quad x\in \mathbb{R}
\label{eq:UBTW-int6}
\end{equation}
in a single space.
The minimal speed of front-type progressive solutions is known exactly for $f(u)$ of the monostable type, $f(u) = au(1-u/K)$, and of the bistable type, $f(u) = u(1-u)(u-a)$.
In particular, the minimal speed for the monostable type is $c_{*}=2\sqrt{f'(0)}=2\sqrt{a}$.
See, for instance, \cite{ito1} and their references.

%非有界進行波について先行研究列挙しつつ
Assume that $f(u)$ is monostable-type nonlinear term.
For instance, let $f(u) = au(1 - u/K)$.
In this case, $f(u) = au$ as $K \to \infty$.
As can be seen from \cite{ito1}, in the reaction-diffusion equation
\begin{equation}
u_{t}=u_{xx}+au, \quad t>0, \quad x\in \mathbb{R},
\label{eq:UBTW-int7}
\end{equation}
we introduce the traveling wave coordinate $u(t,x)=\phi(\xi)$, where $\xi=x-ct$.
Thus, there exists an unbounded traveling wave solution for \eqref{eq:UBTW-int7}.
That is, it can be seen that there exists a monotonically decreasing profile $\phi$ such that 
\begin{equation}
\phi(-\infty)=+\infty, \quad \phi(+\infty)=0.
\label{eq:UBTW-int8}
\end{equation}
Furthermore, it has been shown that the minimal speed is $c_{*}=2\sqrt{a}$, which is the same as the speed of the front-type solution.
In \cite{ito1}, this unbounded traveling wave solution is defined as a traveling wave connecting $\infty$ and $0$.
In addition, for the one-dimensional reaction-diffusion equation \eqref{eq:UBTW-int6} in a single space, assuming that the reaction term $f(u)$ is a $C^{1}$-function satisfying
\begin{equation}
f(0)=0, \quad f'(0)\neq 0, \quad \limsup_{u\to \infty} \dfrac{|f(u)|}{u}<+\infty,
\label{eq:UBTW-int9}
\end{equation}
Ito-Ninomiya explicitly determine the relationship between the existence of traveling wave solutions (unbounded traveling wave solutions) connecting $\infty$ and $0$ and the minimal speed, and they explicitly specify the minimal speed.

%本論文で明らかにしたいこと
In \eqref{eq:UBTW-int1}, $f(u)=u^{3}/(1+su^{2})$ and $f'(0)=0$.
As mentioned above, since \eqref{eq:UBTW-int1} does not satisfy condition \eqref{eq:UBTW-int9} regarding the reaction term of the reaction-diffusion equation considered in \cite{ito1}, it is necessary to investigate it using a different approach.

%戦略の概要：ポアンカレ型コンパクト化
In this paper, we consider the dynamical system on the phase plane induced by the two-dimensional system of ordinary differential equations:
\begin{equation}
\begin{cases}
\phi' = \psi, \\
\psi' = -c\psi - \phi^{3}/(1+s\phi^{2}),
\end{cases}
\quad (\, '=d/d\xi),
\label{eq:UBTW-int10}
\end{equation}
which is equivalent to the second-order equation \eqref{eq:UBTW-int3} obtained by introducing the traveling wave coordinates \eqref{eq:UBTW-int2} for \eqref{eq:UBTW-int1}.
In particular, we consider the dynamical system on $\mathbb{R} \, \cup\, \{\|(\phi, \psi)\|=+\infty\}$ known as the ``dynamics including infinity''.
Since the unbounded traveling wave solutions are defined by their behavior as $\phi \to +\infty$, the problem reduces to investigating the full dynamics including infinity in order to correctly extract this information.
In this paper, we investigate these dynamics including infinity using Poincar\'e-type compactification.
See \cite{FAL, QTW, BIRD, FKPM, UPKPP, Matsue1, Matsue2} or Subsection \ref{sub:UBTW-dyn2} for details on Poincar\'e-type compactification.
This method has been used in the analysis of the Li\'enard equation (see \cite{FAL} and its references) and in the reconstruction of blow-up solutions of ordinary differential equations (ODEs) using dynamical systems theory (see \cite{Matsue1, Matsue2}).
In recent years, the author has applied this method to research aimed at classifying traveling waves in reaction-diffusion equations with nonlinear diffusion. 
For instance, see \cite{FKPM, UPKPP}.

%本論文の結論の概要
In this paper, by investigating the dynamics at infinity revealed by Poincar\'e-type compactification for the traveling wave system \eqref{eq:UBTW-int10}, we establish the existence of unbounded traveling wave solutions connecting $\infty$ and $0$, along with an explicit formula for their minimal speed.
We also provide insight into the existence and profiles of sign-changing unbounded traveling wave solutions.
Furthermore, we show that the information determining the minimal speed does not come from a stable equilibrium at the leading edge of the unbounded traveling wave. 
Instead, it originates from the vicinity of an unstable equilibrium at infinity. 
This differs from conventional linear predictions.
Regarding the minimal speed of unbounded traveling wave solutions, it suggests that information at infinity plays a more significant role than stable information at 0 for asymptotically linear reaction terms.
This analysis provides a strong motivation to investigate other reaction terms and contributes to a new understanding of the existence and characterization of unbounded traveling wave solutions.

%論文の残りの構成
The paper is organized as follows.
The next section presents the main results.
Section \ref{sec:UBTW-dyn} discusses the dynamics including infinity in \eqref{eq:UBTW-int10}.
This is also known as a Poincar\'e disk, which is obtained via Poincar\'e compactification.
In Section \ref{sec:UBTW-pro}, we prove the main result based on the findings from the previous section.
Finally, in Section \ref{sec:UBTW-di}, we summarize the authors' previous findings on unbounded traveling waves, including those presented in this paper. 
We also conclude with various conjectures and observations aimed at future research.

\section{Main results}
\label{sec:UBTW-mr}
The symbol $f(\xi)\sim g(\xi)$ as $\xi\to a$ implies that 
\[
\lim_{\xi\to a}\left| \dfrac{f(\xi)}{g(\xi)} \right|=1.
\]
Before discussing the main results of this paper, we will once again define the traveling wave solution connecting $\infty$ and $0$ (i.e., the unbounded traveling wave solution) (see also \cite{ito1}).

\begin{defn}
\label{def:UBTW1}
A monotonically decreasing, positive profile $\phi=\phi(\xi)$ on $\mathbb{R}$ that satisfies \eqref{eq:UBTW-int3} and \eqref{eq:UBTW-int8} is called an unbounded traveling wave solution of \eqref{eq:UBTW-int1}, or a traveling wave solution connecting $\infty$ and $0$.
\end{defn}

The main results of this paper are as follows.

\begin{thm}
\label{thm:UBTW1}
Let $s > 0$.
Then, for the equation \eqref{eq:UBTW-int1}, there exists a minimal speed (critical value) $0 < c_{*} \in \mathbb{R}$, such that the following holds:
\begin{enumerate}
\item[(i)] If $0 < c < c_{*}$, then there are no positive, unbounded traveling wave solutions.
\item[(ii)] If $c = c_{*}$, then there exists exactly one type of positive, unbounded traveling wave solution with speed $c = c_{*}$ connecting $\infty$ and $0$.
\item[(iii)] If $c > c_{*}$, then there exists exactly two types of positive, unbounded positive traveling wave solutions with speed $c$ connecting $\infty$ and $0$.
\end{enumerate}
Furthermore, the explicit formula for the minimal speed is $c_{*}=2/\sqrt{s}$.
\end{thm}

\begin{thm}
\label{thm:UBTW2}
Let $s > 0$.
Then, for the equation \eqref{eq:UBTW-int1}, there exists a critical value $0 < c_{*} \in \mathbb{R}$, such that the following holds:
\begin{enumerate}
\item[(i)] If $0 < c < c_{*}$, then there exists an unbounded traveling wave solution that oscillates infinitely many times.
\item[(ii)] If $c \ge c_{*}$, then there exists a family of unbounded traveling wave solutions that satisfy the following:
\begin{enumerate}
\item[(ii-1)] The profile $\phi(\xi)$ satisfies \eqref{eq:UBTW-int3} and \eqref{eq:UBTW-int8} on $\mathbb{R}$.
\item[(ii-2)] There exists a unique point $\xi^{*}\in (-\infty, \infty)$ such that $\phi'(\xi)<0$ in $\xi\in (-\infty, \xi^{*})$ and $\phi(\xi^{*})=0$, and $\phi'(\xi)>0$ in $\xi\in (\xi^{*}, +\infty)$.
\end{enumerate}
\end{enumerate}
Furthermore, the explicit formula for $c_{*}$ is $c_{*}=2/\sqrt{s}$.
\end{thm}

\section{Dynamics at finite and infinity}
\label{sec:UBTW-dyn}
In this section, we study the equilibria of \eqref{eq:UBTW-int10} (hereinafter referred to as finite equilibria) and their local stability.
We also investigate the dynamics of $\mathbb{R}^{2}\, \cup\, \{(\phi, \psi) \mid \|(\phi,\psi)\|=+\infty\}$ corresponding to the dynamical system that including infinity.
This can be understood via a Poincar\'e-type compactification.

In \eqref{eq:UBTW-int10}, we impose the following the time-rescale desingularization:
\begin{equation}
d\tilde{s}/d\xi = (1+s \phi^{2})^{-1}.
\label{eq:UBTW-dyn1}
\end{equation}
Therefore, we have 
\begin{equation}
\begin{cases}
\phi_{\tilde{s}} = \psi + s\phi^{2}\psi, \\
\psi_{\tilde{s}} = -c\psi - cs\phi^{2}\psi - \phi^{3}
\end{cases}
\label{eq:UBTW-dyn2}
\end{equation}
with $\phi_{\tilde{s}}=d\phi/d\tilde{s}$ and $\psi_{\tilde{s}}=d\psi/d\tilde{s}$.
This transformation \eqref{eq:UBTW-dyn1} corresponds to multiplying the vector field by $1+s \phi^2$.
The solution curves of the system are invariant under both the parameter $\xi$ and $\tilde{s}$.
See \cite{CK} and the references therein for details.
Therefore, we will examine the dynamics of system \eqref{eq:UBTW-dyn2} at finite equilibria and at infinity, as these are easier to handle analytically.
Equation \eqref{eq:UBTW-dyn2} has the following symmetry:
\begin{equation}
\phi \mapsto -\phi, \quad \psi \mapsto -\psi.
\label{eq:UBTW-dyn3}
\end{equation}

\subsection{Dynamics near finite equilibria}
\label{sub:UBTW-dyn1}
The system \eqref{eq:UBTW-dyn2} has the following equilibrium:
\[
E_{0}: (\phi, \psi)=(0,0).
\]
The Jacobian matrix $J_{0}$ of the vector field \eqref{eq:UBTW-dyn2} at$E_{0}$ is 
\[
J_{0}= \left(\begin{array}{cc}
0 & 1 \\ 0 & -c
\end{array} \right).
\]
The eigenvalues of $J_{0}$ are $0$ and $-c$ and the corresponding eigenvectors are 
\[
{\mathbf{v}}_{1}=(1,0)^{T}, \quad 
{\mathbf{v}}_{2}=(1, -c)^{T},
\]
respectively.
Here, $T$ is the symbol for transpose.
Since the finite equilibrium $E_{0}$ is not hyperbolic, we can understand the center manifold theory to study the dynamics near $E_{0}$.
See \cite{carr, cDNPE, UPKPP} and references therein for details about the center manifold theory. 
We set a matrix $P$ as $P=(\mathbf{v}_{1},\mathbf{v}_{2})$.
Then we obtain 
\begin{align*}
\left(\begin{array}{cc}
\phi_{\tilde{s}} \\
\psi_{\tilde{s}}
\end{array}
\right) &= \left(\begin{array}{cc}
0 & 1 \\ 0 & -c
\end{array}
\right)\left(\begin{array}{cc}
\phi \\
\psi
\end{array}
\right)+\left(\begin{array}{cc}
s\phi^{2}\psi \\
-cs\phi^{2}\psi -\phi^{3}
\end{array}
\right) \\
&= P\left(\begin{array}{cc}
0 & 0 \\
0 & -c
\end{array}
\right)P^{-1}\left(\begin{array}{cc}
\phi \\
\psi
\end{array}
\right)+\left(\begin{array}{cc}
s\phi^{2}\psi \\
-cs\phi^{2}\psi -\phi^{3}
\end{array}
\right).
\end{align*}
Let 
$\left(\begin{array}{cc}
\tilde{\phi} \\
\tilde{\psi}
\end{array}
\right)=P^{-1}\left(\begin{array}{cc}
\phi \\
\psi
\end{array}
\right)$.
We then obtain the following system:
\begin{equation}
\begin{cases}
\tilde{\phi}_{\tilde{s}}=--c^{-1}\tilde{\phi}^{3}-3c^{-1}\tilde{\phi}^{2}\tilde{\psi}-3c^{-1}\tilde{\phi}\tilde{\psi}^{2}-c^{-1}\tilde{\psi}^{3}, \\
\tilde{\psi}_{\tilde{s}}= -c\tilde{\psi}+c^{-1}\tilde{\phi}^{3}+(c^{-1}-2cs)\tilde{\psi}^{3}+(3c^{-1}-cs)\tilde{\phi}^{2}\tilde{\psi}+(3c^{-1}-cs)\tilde{\phi}\tilde{\psi}^{2}.
\end{cases}
\label{eq:UBTW-dyn4}
\end{equation}
The center manifold theory is applicable to study the dynamics of system \eqref{eq:UBTW-dyn4}.
It implies that there exists a function $h(\tilde{\phi})$ satisfying
\[
h(0)=\dfrac{dh}{d\tilde{\phi}}(0)=0
\]
such that the center manifold of the origin for system \eqref{eq:UBTW-dyn4} is locally represented as $\{(\tilde{\phi},\tilde{\psi}) \,|\, \tilde{\psi}(\tilde{s})=h(\tilde{\phi}(s))\}$.
Differentiating it with respect to $\tilde{s}$, we obtain that the approximation of the (graph of) center manifold is
\begin{equation}
\left\{ (\tilde{\phi}, \tilde{\psi}) \,|\, \tilde{\psi}=c^{-2}\tilde{\phi}^{3}+O(\tilde{\phi}^{4}) \right\}.
\label{eq:UBTW-dyn5} 
\end{equation}
Therefore, the dynamics of \eqref{eq:UBTW-dyn4} near $E_0$ is topologically conjugate to the dynamics of the following equation:
\begin{equation}
\tilde{\phi}_{\tilde{s}}(\tilde{s})= -c^{-1} \tilde{\phi}^{3}+O(\tilde{\phi}^{4}).
\label{eq:UBTW-dyn6}
\end{equation}
Hence, we conclude that the approximation of the (graph of) center manifold near $E_{0}$ is
\begin{equation}
\left\{ (\phi, \psi) \mid \psi(\tilde{s})= -\dfrac{1}{c}\left( \phi+\dfrac{1}{c}\psi \right)^{3} +O((\phi+c^{-1}\psi)^{4}) \right\}
\label{eq:UBTW-dyn7}
\end{equation}
and the dynamics of \eqref{eq:UBTW-dyn2} near $E_{0}$ is topologically conjugate to the dynamics of the following equation:
\begin{equation}
\phi_{\tilde{s}}(\tilde{s})=-c^{-1}\phi^{3}+O(\phi^{4}).
\label{eq:UBTW-dyn8}
\end{equation}

\subsection{Short review of Poincar\'e-type compactifications inducing dynamics at infinity}
\label{sub:UBTW-dyn2}
Next, let us extract the dynamics of \eqref{eq:UBTW-dyn2} at infinity.
The system \eqref{eq:UBTW-dyn2} is an asymptotically quasi-homogeneous vector field of type $(1,1)$ and order $3$ at infinity.
See \cite{Matsue1, Matsue2} for the definition and details.
This notion guarantees that the structure of a dynamical system at infinity can be correctly extracted.
Matsue (\cite{Matsue1, Matsue2}) has shown that the scale invariance known as the homogeneity of the vector field must hold at infinity.

Full dynamics $\Phi=\mathbb{R}^{2}\, \cup\, \{(\phi,\psi) \mid \|(\phi,\psi)\|=+\infty \}$ of \eqref{eq:UBTW-dyn2}, including those at infinity, are induced by Poincar\'e-type compactification.
$\Phi$ is hereinafter referred to as a Poincar\'e-type disk.
The outline of this method is described below.
See \cite{FAL, QTW, BIRD, FKPM, UPKPP, Matsue1, Matsue2} and the references therein for details.

Consider a plane $(y_{1},y_{2},y_{3})=(\phi,\psi,1)$ and a sphere in $\mathbb{R}^{3}$, along with their partitions ($k=1,2,3$):
\begin{align*}
\mathbb{S}^{2} &= \{ y \in \mathbb{R}^{3} \, |\, y_{1}^{2} + y_{2}^{2}+y_{3}^{2}=1\}, \quad
\mathbb{S}^{1}=\{y \in \mathbb{S}^{2}\, | \, y_{3}=0\},\\
H_{+} &= \{ y \in \mathbb{S}^{2}\,|\,y_{3}>0\}, \quad
H_{-} = \{ y \in \mathbb{S}^{2}\,|\,y_{3}<0\}, \\
U_{k} &= \{y \in \mathbb{S}^{2} \, | \, y_{k}>0\}, \quad
V_{k} = \{y \in \mathbb{S}^{2} \, | \, y_{k}<0\}.
\end{align*}
Let $\Delta(\phi,\psi)$ be defined as $\Delta(\phi,\psi) = \sqrt{\phi^{2}+\psi^{2}+1}$.
Furthermore, we define the map $f^{\pm}$ from $(\phi, \psi)\in \mathbb{R}^{2}$ to $H^{\pm}$ as 
\[
f^{\pm}(\phi,\psi):= \pm \left( \dfrac{\phi}{\Delta(\phi,\psi)},\dfrac{\psi}{\Delta(\phi,\psi)},\dfrac{1}{\Delta(\phi,\psi)} \right)
\]
and the projection $g^{\pm}_{k}$:
\[
g^{+}_{k} : U_{k} \to \mathbb{R}^{2} \quad {\rm and} 
\quad g^{-}_{k} : V_{k} \to \mathbb{R}^{2} 
\]
as $g^{+}_{k}(y_{1},y_{2},y_{3}) = - g^{-}_{k}(y_{1},y_{2},y_{3})
 = (y_{m}/y_{k}, y_{n}/y_{k})$ for $m<n$ and $m,n\neq k$.
The projected vector fields are obtained as the vector fields on the planes
 \begin{align*}
\overline{U}_{k} = \{y \in \mathbb{R}^{3} \, | \, y_{k} = 1\}, \\
\overline{V}_{k} = \{y \in \mathbb{R}^{3} \, | \, y_{k} = -1\}
\end{align*}
for each local charts $U_{k}$ and $V_{k}$ ($k=1,2,3$).
In addition, we define $(\lambda_{2}, \lambda_{1}) := g^{\pm}_{k}(y)$.

\begin{figure}[t]
\begin{center}
\includegraphics[scale=0.3]{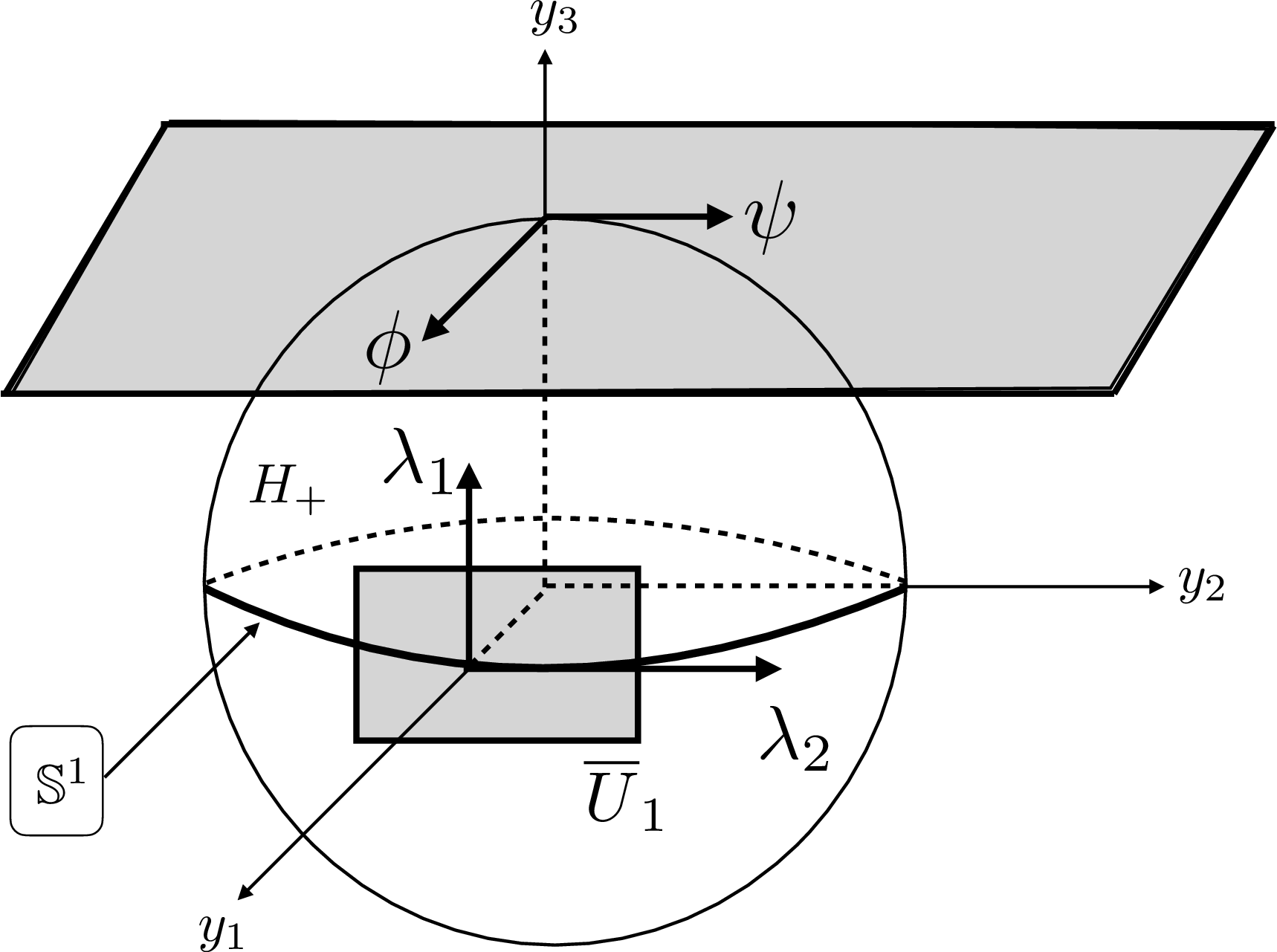}
\caption{Schematic pictures of the locations of the chart $\overline{U}_1$.}
\label{fig:UBTW-dyn0}
\end{center}
\end{figure}

For instance, we consider the case $k=1$.
It follows that the projection from $\mathbb{R}^2$ to $\overline{U}_1$ as 
\[
(g^{+}_{1} \circ f^{+})(\phi,\psi) 
= \left ( \dfrac{\psi}{\phi},\dfrac{1}{\phi}\right) = (\lambda_{2}, \lambda_{1}).
\]
See also Figure \ref{fig:UBTW-dyn0}.
Therefore, we can obtain the dynamics on the local chart $\overline{U}_{1}$ by the change of variables $\phi=1/\lambda_{1}$ and $\psi=\lambda_{2}/\lambda_{1}$.
Throughout this paper, we follow the notations used here for the Poincar\'e-type compactification.

\subsection{Dynamics on the local charts $\overline{U}_{1}$ and $\overline{V}_{1}$}
\label{sub:UBTW-dyn3}
To obtain the dynamics on the chart $\overline{U}_{1}$, we introduce the coordinates $(\lambda_{1}, \lambda_{2})$ given by 
\begin{equation}
\phi=1/\lambda_1, \quad \psi=\lambda_{2}/\lambda_{1}.
\label{eq:UBTW-dyn9}
\end{equation}
Then 
\begin{equation}
\begin{cases}
d\lambda_{1}/d\tilde{s}
= -\lambda_{1}\lambda_{2} - s\lambda_{1}^{-1}\lambda_{2},
\\
d\lambda_{2}/d\tilde{s}
= -c\lambda_{2}-cs \lambda_{1}^{-2}\lambda_{2} - \lambda_{1}^{-2}-\lambda_{2}^{2}-s\lambda_{1}^{-2}\lambda_{2}^{2}
\end{cases}
\label{eq:UBTW-dyn10}
\end{equation}
holds.
By using the time-scale desingularization $d\tau/d\tilde{s}=\lambda_{1}^{-2}$, we obtain
\begin{equation}
\begin{cases}
d\lambda_{1}/d\tau 
= -\lambda_{1}^{3}\lambda_{2} - s\lambda_{1}\lambda_{2}, \\
d\lambda_{2}/d\tau
= -c\lambda_{1}^{2}\lambda_{2} -cs \lambda_{2} -1-\lambda_{1}^{2}\lambda_{2}^{2} - s\lambda_{2}^{2}.
\end{cases}
\label{eq:UBTW-dyn11}
\end{equation}
The equilibria of \eqref{eq:UBTW-dyn11} on $\{\lambda_{1}=0\}$, which correspond to the equilibria at infinity, can be classified as follows:
\begin{enumerate}
\item[(i)] 
When $c^{2}s^{2}-4s>0$ (i.e., $c>2/\sqrt{s}$), the equilibria are
\[
E_{1}: (\lambda_{1}, \lambda_{2})=(0, M_{-}), \quad
E_{2}: (\lambda_{1}, \lambda_{2})=(0, M_{+}), \quad 
M_{\pm}=\dfrac{-cs\pm \sqrt{c^{2}s^{2}-4s}}{2s}<0.
\]
\item[(ii)] 
When $c^{2}s^{2}-4s=0$ (i.e., $c=2/\sqrt{s}$), the system \eqref{eq:UBTW-dyn11} has the equilibrium $E_{3}: (\lambda_{1}, \lambda_{2})=(0, -s^{-1/2})$.
\item[(iii)] 
When $c^{2}s^{2}-4s<0$ (i.e., $0<c<2/\sqrt{s}$), the system \eqref{eq:UBTW-dyn11} on $\{\lambda_{1}=0\}$ has no equilibria.
\end{enumerate}
The linearized matrices $J_{i}$ ($i=1,2,3$) of the vector field \eqref{eq:UBTW-dyn11} at each equilibrium $E_{i}$ are 
\[
J_{1}=\left(\begin{array}{cc}
-M_{-} & 0 \\ 0 & \sqrt{c^{2}s^{2}-4s}
\end{array}\right), \quad
J_{2}=\left(\begin{array}{cc}
-M_{+} & 0 \\ 0 & -\sqrt{c^{2}s^{2}-4s}
\end{array}\right),\quad
J_{3}=\left(\begin{array}{cc}
\sqrt{s} & 0 \\ 0 & 0
\end{array}\right).
\]
Hence, we can conclude that the equilibrium $E_{1}$ is a source and $E_{2}$ is a saddle.
Since the matrix $J_{3}$ has the zero eigenvalue, $E_{3}$ is not hyperbolic.
As in Subsection \ref{sub:UBTW-dyn1}, applying the center manifold theorem to the dynamical system near $E_3$ yields
\begin{equation}
\{(\lambda_{1}, \lambda_{2}) \mid \lambda_{1}=0\}
\label{eq:UBTW-dyn12}
\end{equation}
as an approximation of the center manifold.
In addition, the dynamics of \eqref{eq:UBTW-dyn11} near $E_3$ is topologically conjugate to the dynamics of the following equation:
\begin{equation}
\dfrac{d\lambda_{2}}{d\tau}= -s\lambda_{2}^{2}-cs\lambda_{2}-1.
\label{eq:UBTW-dyn13}
\end{equation}

Similarly, the dynamics on the local chart $\overline{V}_{1}$ can be obtained by introducing the transformation 
\begin{equation}
\phi= -1/\lambda_1, \quad \psi= -\lambda_{2}/\lambda_{1}.
\label{eq:UBTW-dyn14}
\end{equation}
This transformation yields \eqref{eq:UBTW-dyn9}.
It arises from the symmetry \eqref{eq:UBTW-dyn3} of \eqref{eq:UBTW-dyn2}.

\subsection{Dynamics on the local charts $\overline{U}_{2}$ and $\overline{V}_{2}$}
\label{sub:UBTW-dyn4}
Since the system \eqref{eq:UBTW-dyn2} has the symmetry \eqref{eq:UBTW-dyn3}, it is sufficient to consider only the case of the local chart $\overline{U}_{2}$.
The change of coordinates 
\begin{equation}
\phi=\lambda_{2}/\lambda_{1}, \quad
\psi=1/\lambda_{1}
\label{eq:UBTW-dyn15}
\end{equation}
and time-rescaling $d\tau/d\tilde{s}=\lambda_{1}^{-2}$ to study the dynamics on the local chart $\overline{U}_{2}$ yields 
\begin{equation}
\begin{cases}
d\lambda_{1}/d\tau
=c\lambda_{1}^{3}+cs\lambda_{1}\lambda_{2}^{2}+\lambda_{1}\lambda_{2}^{3},
\\
d\lambda_{2}/d\tau
=\lambda_{1}^{2}+s\lambda_{2}^{2}+c\lambda_{1}^{2}\lambda_{2}+cs\lambda_{2}^{3}+\lambda_{2}^{4}.
\end{cases}
\label{eq:UBTW-dyn16}
\end{equation}
The equilibria on $\{\lambda_{1}=0\}$ are classified as follows:
\begin{enumerate}
\item[(i)] 
When $c^{2}s^{2}-4s>0$ (i.e., $c>2/\sqrt{s}$), the system \eqref{eq:UBTW-dyn16} has the following three types of equilibria:
\[
E_{4}: (\lambda_{1}, \lambda_{2})=(0, sM_{-}), \quad
E_{5}: (\lambda_{1}, \lambda_{2})=(0, sM_{+}), \quad 
E_{6}: (\lambda_{1}, \lambda_{2})=(0,0).
\]
\item[(ii)] 
When $c^{2}s^{2}-4s=0$ (i.e., $c=2/\sqrt{s}$), this system has an equilibrium $E_{7}: (\lambda_{1}, \lambda_{2})=(0, -\sqrt{s})$ in addition to $E_{6}$.
\item[(iii)] 
When $c^{2}s^{2}-4s<0$ (i.e., $0<c<2/\sqrt{s}$), only the equilibrium $E_{6}$ exists.
\end{enumerate}
As will be discussed in the next section, since the equilibria $E_{4,5,7}$ coincide with the equilibria on $\overline{V}_{1}$ in the $\phi\psi$-plane where $\|(\phi, \psi)\| = +\infty$.
Therefore, the conclusions regarding their local stability are the same as those in Subsection \ref{sub:UBTW-dyn3}.

In the remainder of this section, we only need to study the dynamical system near $E_{6}$.
Since the linearized matrix of the system \eqref{eq:UBTW-dyn16} in $E_{6}$ has the double zero eigenvalues, the dynamics of the system near $E_{6}$ can be understood by studying the dynamics in the blow-up vector field obtained by introducing the blow-up coordinates: 
\begin{equation}
\lambda_{1} = r \bar{\lambda}_{1}, \quad 
\lambda_{2} = r \bar{\lambda}_{2}.
\label{eq:UBTW-dyn17}
\end{equation}
See \cite{FAL, CK} and the references therein for details on the directional blow-up technique.

Geometrically, we only need to study the dynamics on $\{r=0\}$ that correspond to the dynamics of $E_6$.
Therefore, since we are focusing on $\lambda_{1} \ge 0$, we will investigate the dynamics on $\{r = 0\}$ in the three local coordinates: $\bar{\lambda}_1 = 1$ and $\bar{\lambda}_2 = \pm 1$.
By integrating these three dynamics, we can obtian the dynamics restricted to $\lambda_1 \ge 0$ on $E_6$.

\subsubsection{Dynamics on the chart $\{\bar{\lambda}_{1}=1 \}$}
\label{sub:UBTW-dyn4-1}
By the change of coordinates $\lambda_{1}=r$, $\lambda_{2}=r\bar{\lambda}_{2}$ and the time-rescaling $d\eta/d\tau=r$, 
\begin{equation}
\begin{cases}
dr/d\eta
= cr^{2}+cs r^{2}\bar{\lambda}_{2}^{2}+r^{3}\bar{\lambda}_{2}^{3}, 
\\
d\bar{\lambda}_{2}/d\eta
=1+s\bar{\lambda}_{2}^{2}
\end{cases}
\label{eq:UBTW-dyn18}
\end{equation}
holds.
There is no equilibrium, and $d\lambda_{2}/d\tau > 0$ holds on $\{r=0\}$.

\subsubsection{Dynamics on the chart $\{\bar{\lambda}_{2}=\pm1 \}$}
\label{sub:UBTW-dyn4-2}
Using the change of coordinates $\lambda_{1}=r\bar{\lambda}_{1}$, $\lambda_{2}=r$ and time-rescaling $d\eta/d\tau=r$, the dynamics on the local chart $\{\bar{\lambda}_{2}=1 \}$ yields
\begin{equation}
\begin{cases}
dr/d\eta
= r\bar{\lambda}_{1}^{2}+sr+cr^{2}\bar{\lambda}_{1}^{2}+csr^{2}+r^{3}, 
\\
d\bar{\lambda}_{1}/d\eta
=-\bar{\lambda}_{1}^{3}-s\bar{\lambda}_{1}.
\end{cases}
\label{eq:UBTW-dyn19}
\end{equation}
The system \eqref{eq:UBTW-dyn19} has the equilibrium $(r, \bar{\lambda}_{1})=(0,0)$.
The eigenvalues of the linearized matrix at the origin (equilibrium) are $s$ and $-s$.
Hence, the origin is a saddle.

Similarly, the change of coordinates $\lambda_{1}=r\bar{\lambda}_{1}$ and $\lambda_{2}=r$ yields
\begin{equation}
\begin{cases}
dr/d\eta
= -r\bar{\lambda}_{1}^{2}-sr+cr^{2}\bar{\lambda}_{1}^{2}+csr^{2}-r^{3}, 
\\
d\bar{\lambda}_{1}/d\eta
=\bar{\lambda}_{1}^{3}+s\bar{\lambda}_{1}.
\end{cases}
\label{eq:UBTW-dyn20}
\end{equation}
It is the dynamics on the local chart $\{\bar{\lambda}_{2}=-1 \}$ with the time-rescaling$d\eta/d\tau=r$.
The only equilibrium is the origin, and by the same discussion as above, we can see that the origin is a saddle.

\begin{figure}[t]
\begin{center}
\includegraphics[scale=0.35]{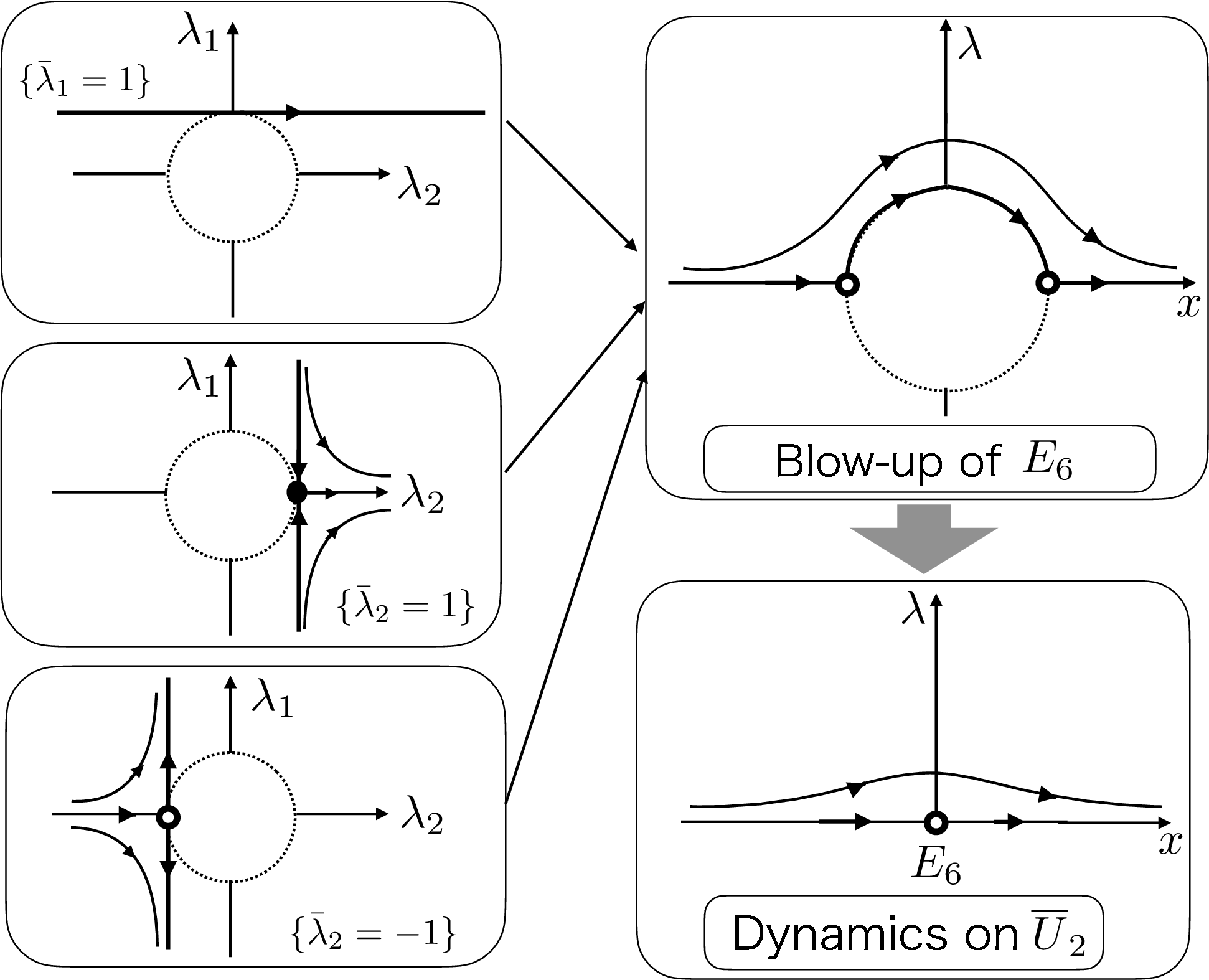}
\caption{Schematic pictures of the dynamics of the blow-up vector fields and dynamics near $E_{6}$.}
\label{fig:UBTW-dyn1}
\end{center}
\end{figure}

Therefore, the dynamics of $E_{6}$ restricted to $\lambda_{1} \ge 0$ are as shown in Figure \ref{fig:UBTW-dyn1}.

\subsection{Dynamics and connecting orbits on the Poincar\'e-type disk}
\label{sub:UBTW-dyn6}
\begin{figure}[t]
\begin{center}
\includegraphics[scale=0.37]{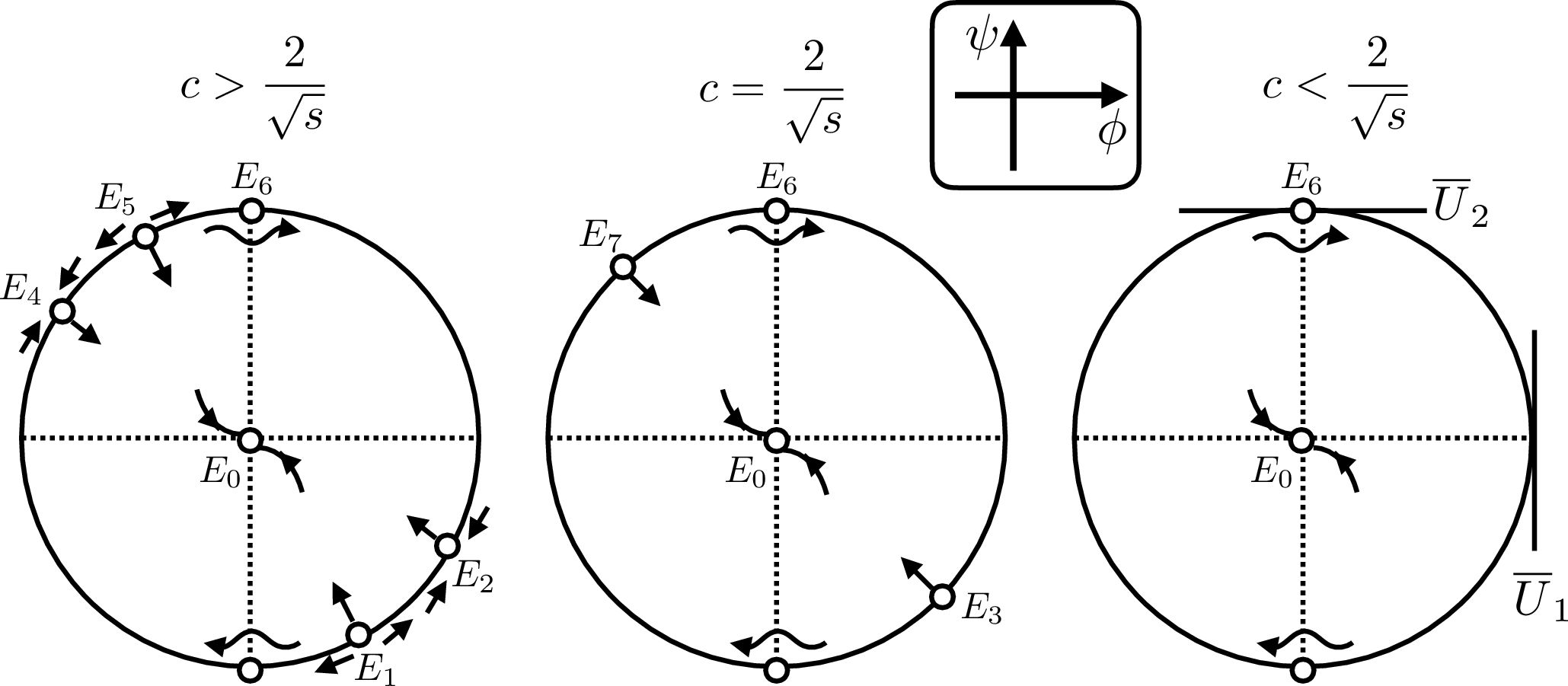}
\caption{Schematic pictures of the dynamics on the Poincar\'e-type disk $\Phi$ for \eqref{eq:UBTW-dyn2}.}
\label{fig:UBTW-dyn2}
\end{center}
\end{figure}

The dynamical system for \eqref{eq:UBTW-dyn2} on $\Phi$ is shown in Figure \ref{fig:UBTW-dyn2}.
It combines information from the dynamical system near the finite equilibria from Subsection \ref{sub:UBTW-dyn1} with the dynamical system at infinity from Subsections \ref{sub:UBTW-dyn3} and \ref{sub:UBTW-dyn4}.
The goal of this section is to prove the existence of connecting orbits between equilibria on $\Phi$ and to classify them.
These orbits correspond to traveling wave solutions (including the unbounded traveling wave solutions) of the original partial differential equation \eqref{eq:UBTW-int1}.

Let us prepare the symbols used in this paper as follows: 
\begin{itemize}
\item 
$\mathcal{W}^{cs}(E_{0})$ denotes the center stable manifold of $E_{0}$ in the dynamical system \eqref{eq:UBTW-dyn2}.
\item 
$\mathcal{W}^{cu}(E_{3})$ denotes the center unstable manifold of $E_{3}$ in the dynamical system \eqref{eq:UBTW-dyn11}.
$E_3$ is an equilibrium on $\{\|(\phi,\psi)\|=+\infty\}$ in \eqref{eq:UBTW-dyn2}, and $\mathcal{W}^{cu}(E_3)$ is also the center unstable manifold in \eqref{eq:UBTW-dyn2}.
\item 
$\mathcal{W}^{u}(E_{1})$ denotes the unstable manifold of $E_{1}$ in the dynamical systems \eqref{eq:UBTW-dyn11} and \eqref{eq:UBTW-dyn2}.
\item 
$\mathcal{W}^{u}(E_{2})$ denotes the unstable manifold of $E_{2}$ in the dynamical systems \eqref{eq:UBTW-dyn11} and \eqref{eq:UBTW-dyn2}.
\end{itemize}

\begin{figure}[t]
\begin{center}
\includegraphics[scale=0.37]{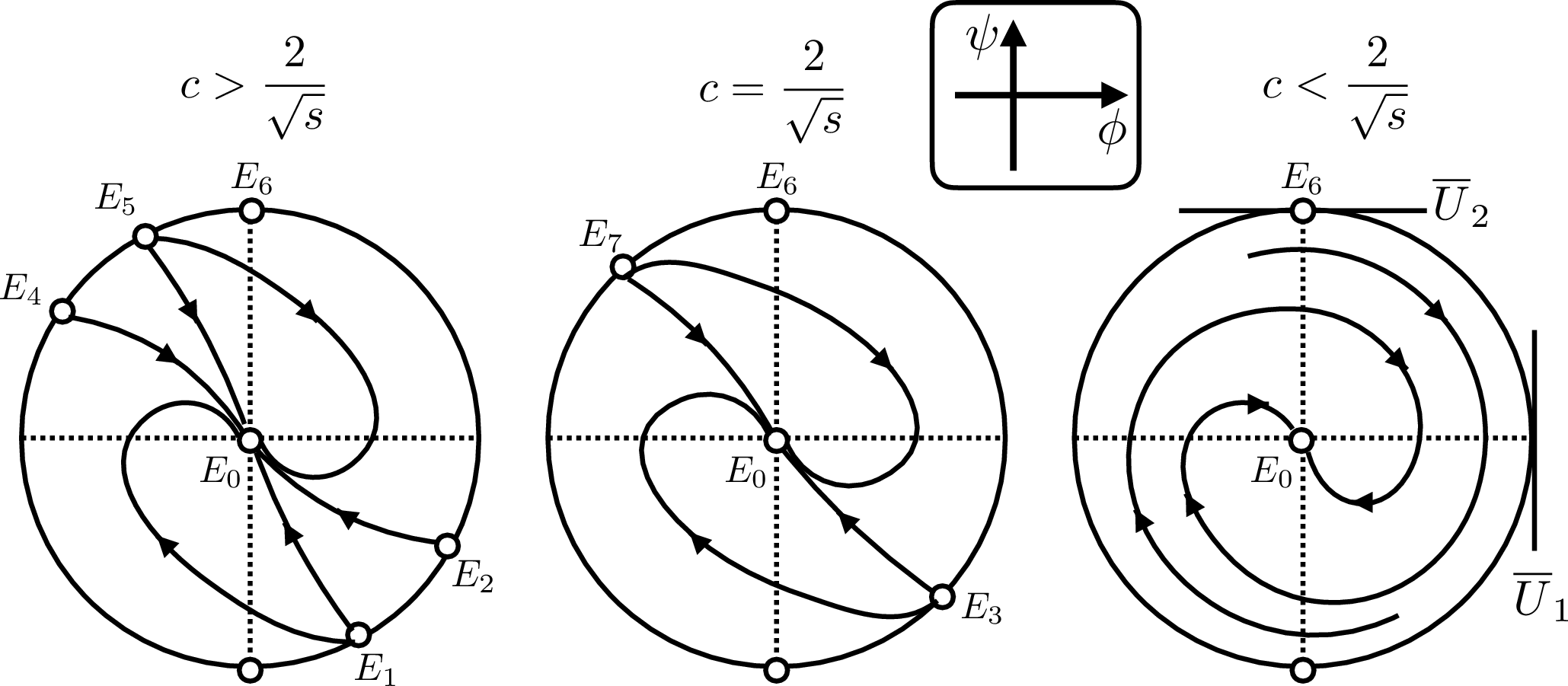}
\caption{Schematic pictures of the dynamics and connecting orbits on the Poincar\'e-type disk $\Phi$ for \eqref{eq:UBTW-dyn2} and \eqref{eq:UBTW-int10}.}
\label{fig:UBTW-dyn3}
\end{center}
\end{figure}

First, consider the case when $c < 2/\sqrt{s}$.
For the dynamical system of \eqref{eq:UBTW-dyn2}, consider the time-reversed system with $\tilde{s} \mapsto -\tilde{s}$.
It then becomes clear that a trajectory starting from an initial point slightly off from $E_{0}$ can only wind its way to infinity.
Consequently, we obtain the diagram on the right of Figure \ref{fig:UBTW-dyn3}.

Second, we consider $c=2/\sqrt{s}$.
Take an initial point on $\mathcal{W}^{cu}(E_3)$.
Since no orbits sink into $E_{6}$, and by the symmetry \eqref{eq:UBTW-dyn3} and the Poincar\'e-Bendixson theorem (see \cite{Wiggins}), we see that orbits starting from this point must go to $\mathcal{W}^{cs}(E_0)$.
Thus, the existence of a connecting orbit between $E_{3}$ and $E_{0}$ is shown.
In $\Phi_{1}=\{(\phi, \psi) \mid \phi=0,\,\,\psi<0\} \subset \Phi$, $\phi_{\tilde{s}}<0$ holds.
In addition, in $\Phi_{2}=\{(\phi, \psi) \mid \phi<0,\,\,\psi<0\} \subset \Phi$, $\phi_{\tilde{s}}<0$ and $\psi_{\tilde{s}}>0$ hold.
Hence, we can see that there also exists a trajectory that starts from a point on $\mathcal{W}^{cu}(E_3)$, passes through the $\psi$-axis from the fourth quadrant of the $\phi\psi$-plane, enters the third quadrant, passes through the $\phi$-axis, and going to $\mathcal{W}^{cs}(E_0)$ in the second quadrant.
See also the central figure in Figure \ref{fig:UBTW-dyn3}.

Finally, consider the case where $c > 2/\sqrt{s}$.
By the same discussion as above, orbits starting from the point of $\mathcal{W}^{u}(E_1)$ and $\mathcal{W}^{u}(E_2)$ can only reach a point on $\mathcal{W}^{cs}(E_0)$.
Using the symmetry \eqref{eq:UBTW-dyn3}, we can prove the existence of connecting orbits between $E_1$ and $E_0$, between $E_2$ and $E_0$, between $E_4$ and $E_0$, and between $E_5$ and $E_0$.
Since we are considering the time-scale transformation \eqref{eq:UBTW-dyn1}, the direction of the connecting orbits remains invariant under the parameters $\xi$ and $\tilde{s}$.
Therefore, it follows that these connecting orbits in \eqref{eq:UBTW-dyn2} are inherited by \eqref{eq:UBTW-int10}.
Combining these results yields Figure \ref{fig:UBTW-dyn3}.

\section{Proof of theorems}
\label{sec:UBTW-pro}
\subsection{Proof of Theorem \ref{thm:UBTW1}}
\label{sub:UBTW-pro1}
Unbounded traveling wave solutions correspond to trajectories that start from equilibria on the $\{ \|(\phi,\psi)\|=+\infty \}$ in the $\phi\psi$ plane on $\Phi$.
From the discussion in Subsection \ref{sub:UBTW-dyn6}, it follows that when $0<c<2/\sqrt{s}=c_{*}$, all orbits pass through the $\psi$-axis where $\phi=0$.
Therefore, there are no positive-valued unbounded traveling wave solutions.

When $c=c_{*}$, the only connecting orbits existing in the region $\phi>0$ within $\Phi$ are those that do not pass through the $\psi$-axis between $E_3$ and $E_0$.
These orbits correspond to positive, unbounded traveling wave solutions.
Furthermore, since this orbit corresponds to $\phi$ decreasing monotonically with respect to $\xi$, it follows that the corresponding positive unbounded traveling wave solution also possesses monotonicity.

When $c > c_* = 2/\sqrt{s}$, the connecting orbits within the region $\phi > 0$ correspond to the family of connecting orbits between $E_1$ and $E_0$ and the connecting orbits between $E_2$ and $E_0$, respectively.
There are positive, unbounded traveling wave solutions.
As mentioned above, the solution $\phi$ is also monotonicity with respect to $\xi$.
Therefore, when $c > c_*$, we see that there exist two types of positive, monotonicity, nonnegative unbounded wave solutions that connect $\infty$ and $0$.
Thus, this completes the proof of Theorem \ref{thm:UBTW1}.
\qed
\\

\subsection{Proof of Theorem \ref{thm:UBTW2}}
\label{sub:UBTW-pro2}
A sign-changing unbounded traveling wave solution corresponds to a solution $\Phi$ that starts from an equilibrium on $\{ \|(\phi,\psi)\|=+\infty \}$ and whose trajectory passes through the $\phi$-axis ($\psi = 0$).

When $0 < c < c_{*}= 2/\sqrt{s}$, all orbits pass through the $\phi$-axis and wrap around infinity.
Thus, it follows that there exists an unbounded traveling wave solution that oscillates infinitely.

When $c = c_*$, there exists a family of orbits between $E_3$ and $E_0$ that crosses the $\phi$-axis exactly once at $\phi < 0$ and changes from $\psi < 0$ to $\psi > 0$.

When $c > c_*$, the family of connecting orbits between $E_1$ and $E_0$ passes through the $\phi$-axis exactly once at $\phi < 0$.
The fact the $\phi$-axis is crossed exactly once implies that there exists a unique point $\xi^{*}$ in $(-\infty, \infty)$ such that $\psi(\xi^{*}) = 0$.
Thus, the claim of Theorem \ref{thm:UBTW2} (ii) is shown.
Thus, this completes the proof of Theorem \ref{thm:UBTW2}.
\qed
\\

\section{Discussion}
\label{sec:UBTW-di}
Equation \eqref{eq:UBTW-int1} has no bounded traveling wave solutions, including front-type ones.
All existing traveling wave solutions are unbounded.
The physical context of the original equation is beyond the scope of this paper and is therefore not addressed here.
Nevertheless, the main results of this paper are important for the model's validity and comparison with phenomena and experiments.

%得られた結果の概要
In this paper, we considered traveling wave solutions to the one-dimensional reaction-diffusion equation \eqref{eq:UBTW-int1}, which features an asymptotically linear reaction term and a linear diffusion term involving a saturation parameter.
By introducing the traveling wave coordinates \eqref{eq:UBTW-int2}, the problem reduces to the second order ordinary differential equation \eqref{eq:UBTW-int3}.
Using the Poincar\'e-type compactification, we have characerized the full dynamics on the $\phi\psi$-phase plane, including the behavior at infinity.
This approach allowed us to classify all trajectories on the $\phi\psi$-plane corresponding to traveling waves.
Our analysis demonstrates that front-type traveling wave solutions do not exist and that all existing traveling wave solutions are unbounded.
Furthermore, by analyzing the linearized equations near the equilibria at infinity, we revealed a relationship between the minimal speed and the existence of monotonic unbounded traveling wave solutions.
These results provide a new perspective that differs from that of \cite{ito1}.

%市田の非線形拡散での非有界進行波との関連
Furthermore, in recent years, the author has presented results regarding the classification of traveling waves for parabolic equations with nonlinear diffusion, as seen in \cite{FKPM, UPKPP}.
These works provide a classification of traveling waves when the reaction term is of the Fisher-KPP type (monostable) and the diffusion term is either of the porous-medium type or of the power-law nonlinear type.

In \cite{UPKPP}, it was proven that for any wave speed $c > 0$, there exists an unbounded traveling wave solution connecting $\infty$ to $1$.
According to \cite{FKPM}, unbounded traveling wave solutions that connect $\infty$ to $1$ exist for any wave speed $c > 0$. 
Additionally, there is a minimal speed $c_*$ for which unbounded traveling wave solutions that connect $\infty$ and $0$ exist.
These findings suggest that nonlinear diffusion leads to a rich variety of profiles for unbounded traveling wave solutions.
Integrating these previous results with the findings of the present paper gives rise to a new open question: how does the profile of an unbounded traveling wave solution change when the linear diffusion in equation \eqref{eq:UBTW-int1} is replaced with nonlinear diffusion?
Investigating this problem is of significant interest.
Thus, the analysis presented in this paper not only builds upon the work of \cite{ito1} but also opens up the possibility of elucidating the relationship between diffusion and unbounded traveling wave solutions by examining the dynamics at infinity in the phase space.

\section*{Acknowledgments}
This work was partially supported by JSPS KAKENHI Grant Number 25K17306.

%%%%%%%%%%%%%%%%%%%%%%%%%%%%%%%%%%%%%%%%%%%%%%%%%%

\end{document}